\documentclass{article}

\usepackage[T1]{fontenc}

\usepackage{graphicx}
\usepackage{subfigure}

\usepackage{amsmath}
\usepackage{amssymb}

\usepackage{url}

\usepackage{mathtools}
\usepackage{tikz}
\usepackage{algorithm}
\usepackage{algorithmicx}
\usepackage{algpseudocode}

\newcommand{\amplitude}{{\alpha}}
\newcommand{\damplitude}{{a}}
\newcommand{\J}{{J}}

\newcommand{\modeIdx}{{k}}
\newcommand{\modes}{{\psi}}
\newcommand{\dmodes}{{w}}
\newcommand{\iter}{{p}}
\newcommand{\nShifts}{{N_\mathrm{s}}}
\newcommand{\rhoref}{{\rho_{\text{ref}}}}
\newcommand{\shift}{{\Delta}}
\newcommand{\dshift}{{d}}
\newcommand{\rVec}{{\boldsymbol{r}}}

\newcommand{\T}{{T}}
\newcommand{\Tinf}{{\mathcal{T}}}
\newcommand{\Tc}{{T_{\mathrm{c}}}}
\newcommand{\Tcinf}{{\mathcal{T}_{\mathrm{c}}}}
\newcommand{\Tper}{{T_{\mathrm{per}}}}
\newcommand{\Tperinf}{{\mathcal{T}_{\mathrm{per}}}}

\newcommand{\eg}{{e.\,g.}}
\newcommand{\ie}{{i.\,e.}}

\begin{document}            
\title{Model Reduction for a Pulsed Detonation Combuster via Shifted Proper Orthogonal Decomposition }
\author{Philipp Schulze\thanks{Institut f\"ur Mathematik, TU Berlin, Germany, \texttt{$\{$pschulze,mehrmann$\}$@math.tu-berlin.de}.} \and Julius Reiss\thanks{Institut f\"ur Str\"omungsmechanik und Technische Akustik, TU Berlin, Germany, \texttt{reiss@tnt.tu-berlin.de}.}
\and Volker Mehrmann\footnotemark[1]}

\maketitle              
\begin{abstract}
We propose a new algorithm to compute a shifted proper orthogonal decomposition (sPOD) for systems dominated by multiple transport velocities.
The sPOD is a recently proposed mode decomposition technique which overcomes the poor performance of classical methods like the proper orthogonal decomposition (POD) for transport-dominated phenomena.
This is achieved by identifying the transport directions and velocities and by shifting the modes in space to track the transports. Our new algorithm carries out a residual minimization in which the main computational cost arises from solving a nonlinear optimization problem scaling with the snapshot dimension.
We apply the algorithm to snapshot data from the simulation of a pulsed detonation combuster and observe that very few sPOD modes are sufficient to obtain a good approximation.
For the same accuracy, the common POD needs ten times as many modes and, in contrast to the sPOD modes, the POD modes do not reflect the moving front profiles properly.
\end{abstract}


\section{Introduction\label{sec:intro}}

Model reduction is an essential requirement in almost all areas of science and technology to obtain efficient multi-parameter simulations and, in particular, optimization and control methods, see e.g.
\cite{Ant05,BenGW15,HesRS16}.
Often the full-order model (FOM) arises from a semi-discretization in space of a partial differential equation (PDE) and the state dimension scales with the number of grid points which is typically large.
However, one  is  usually not interested in a detailed description of the complete dynamics but often only in a low-dimensional manifold where the solution of interest approximately evolves.

Model reduction for nonlinear dynamical systems is often based on mode decomposition techniques as the proper orthogonal decomposition (POD) \cite{BenGW15,BerHL93,Vol01} or the dynamic mode decomposition \cite{KutBBP16,SchS08}.
Standard mode decomposition techniques are based on the concept of representing the unknown solution as a linear combination of modes.
More precisely, let $q$ be a function in space $x$ and time $t$ representing the state of the dynamical system, then a common model reduction ansatz is an approximation
\begin{equation}
	\label{eq:LCmodes}
	q\left(x,t\right) \approx \sum\limits_{k=1}^r \amplitude_k\left(t\right)\modes_k\left(x\right)
\end{equation}
with space-dependent modes $\modes_k$, time-dependent coefficients, or amplitudes, $\amplitude_k$,
and $r$ is the number of modes.

While the amplitudes typically become the unknowns of the reduced-order model, the modes have to be determined in advance. To determine the modes, one typically simulates the system
and computes space- and time-discrete snapshots of a numerical approximation $q_m$ which are stored in a snapshot matrix $X\in\mathbb{R}^{m\times n}$, \ie, $[X]_{ij} = q_m(x_i,t_j) \approx q(x_i,t_j)$ for $i=1,\ldots,m$ and $j=1,\ldots,n$. With the coefficients of the snapshot matrix one obtains a
discrete analogue of  \eqref{eq:LCmodes} as
\begin{equation}
	\label{eq:snapshotDecomposition}
	\left[X\right]_{j} \approx \sum\limits_{k=1}^r \damplitude_{k,j}\dmodes_k
\end{equation}
for $j=1,\ldots,n$, where $\left[X\right]_{j}$ denotes the $j$th column of $X$, $\dmodes_k\in\mathbb{R}^m$ are coefficient vector representations of the modes $\modes_k$, and $\damplitude_{k,j}$ are the corresponding amplitudes at time point $t_j$.

A classical way to obtain modes and amplitudes is the POD which is based on a singular value decomposition (SVD) of the snapshot matrix $X$. The POD  representation is optimal in the sense that it minimizes the residual in the discrete representation \eqref{eq:snapshotDecomposition}. The resulting reduced-order model is obtained as projection onto the span of the so obtained modes.

A common assumption is that POD delivers a good approximation of the form \eqref{eq:LCmodes} or \eqref{eq:snapshotDecomposition} with a small number $r$.
In many applications this assumption is valid and model reduction schemes like POD  lead to models with dimensions that are orders of magnitude smaller than those of the full-order model \cite{HolLB96}.

However, when the dynamics of the system features structures with high gradients that are propagating through the domain, then schemes of the form \eqref{eq:LCmodes} typically need a large number of modes to approximate the dynamics well, and hence model order reduction is not very effective.
To overcome this difficulty, recently, there have been several suggestions for model reduction of such transport-dominated systems.
In \cite{OhlR13} the authors use ideas of symmetry reduction to decompose the solution into a frozen profile and a translation group accounting for the transport.
The advantages over standard model reduction schemes are demonstrated by means of the Burgers' equation.
In \cite{RimML17_ppt} the authors present a method which is able to decompose multiple transport phenomena. The main ingredients are SVDs of several shifted snapshot matrices combined with a greedy algorithm. The method is cheap to apply but it often needs more shifted modes than necessary, as illustrated with results for the linear wave equation.
For further references on model reduction for transport-dominated problems, see \cite{AbgAC16,CagMS16_ppt,GerL14,IolL14,MojB17,SesS16}.
Most of these approaches consider transport-dominated systems with only one transport velocity and assume periodic boundary conditions.
However, in many applications, multiple transport velocities are encountered, \eg, by different waves propagating through the domain.
To deal with such phenomena,  in \cite{ReiSSM17_ppt} the shifted POD (sPOD) method has been proposed  to obtain mode decompositions suitable for multiple transport phenomena. This new technique differs from \eqref{eq:LCmodes} by shifting the modes in space into different reference frames according to the different transports of the system, \ie,
\begin{equation}
\label{eq:shiftedSnapshotDecomposition}
q(x,t)\approx  \sum\limits_{\ell=1}^{\nShifts} \Tperinf\left(\shift^{\ell}(t)\right) \sum\limits_{k=1}^{r^\ell} \amplitude^{\ell}_{\modeIdx}(t)\modes^{\ell}_{\modeIdx}(x)
\end{equation}
where $\Tperinf(\cdot)$ is a shift operator defined on a periodic domain $[0,L]$ via
\begin{equation*}
	\label{eq:periodicShiftOperator}
	\Tperinf\left(\shift(t)\right)f\left(x,t\right)\vcentcolon=f\left(\left(x+\shift (t)\right)\,\text{mod}\, L,t\right),
\end{equation*}
$\nShifts$ denotes the number of shifted reference frames, $\text{mod}$ denotes the \emph{modulo operator} reflecting the periodicity of the domain,
and $\shift^\ell (t)$ are time-dependent shifts which track the locations of, \eg, different wave profiles over time.
Similar as in the POD, one obtains a discrete analogue of \eqref{eq:shiftedSnapshotDecomposition} via
\begin{equation}
\label{eq:shiftedSnapshotDecomposition_discrete}
\left[X\right]_{j} \approx  \sum\limits_{\ell=1}^{\nShifts} \Tper\left(\dshift^{\ell}_{j}\right) \sum\limits_{k=1}^{r^\ell} \damplitude^{\ell}_{k,j}\dmodes^{\ell}_{k}
\end{equation}
for $j=1,\ldots,n$, where $\Tper$ is a discrete approximation of $\Tperinf$ and $\dshift^{\ell}_{j}$ are shifts at discrete time points $t_j$.
In \cite{ReiSSM17_ppt} a heuristic algorithm is proposed to compute a decomposition of the form \eqref{eq:shiftedSnapshotDecomposition_discrete} in an iterative procedure, and it has been demonstrated that this approach is very successful for several examples including two separating vortex pairs and the linear wave equation. In the latter case the method needs less modes than other methods such as \eg\ \cite{RimML17_ppt} and also retrieves the known analytic solution.

In this paper, we propose an optimization procedure to compute an \emph{optimal} decomposition of the form \eqref{eq:shiftedSnapshotDecomposition_discrete}. To this end, we generalize the cost functional
that is used to obtain the optimality of the POD method to the sPOD setting. We first consider
the optimization on the infinite-dimensional level, see Sec. \ref{sec:optimalspod}, and then present an algorithm which computes the decomposition in the fully discrete setting, see Sec.\ \ref{sec:optimization} and \ref{sec:algorithm}.
The computational cost is higher than for the method in \cite{ReiSSM17_ppt} but the obtained approximations are locally optimal in the sense that a  residual is minimized.

The focus of our work is on obtaining an optimal  mode decomposition which then can be used for the  construction of a reduced-order model, \eg, by a Galerkin projection. A rigorous treatment of non-periodic boundary conditions is also discussed elsewhere.

To demonstrate the efficiency of the new approach, we present results for a pulsed detonation combuster (PDC). The snapshot data originate from a data assimilation, cf. \cite{GraLRPSM17}, and exhibit multiple transport phenomena which interact nonlinearly with each other and with the boundary.

\section{Optimal sPOD Approximation}\label{sec:optimalspod}

\begin{figure}[htb]
	\centering
	\includegraphics [width=0.325\linewidth]{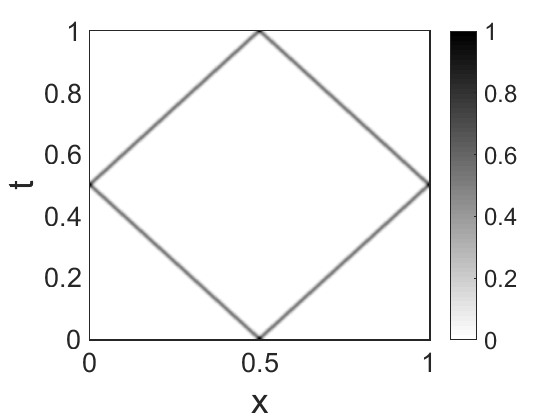}
	\includegraphics [width=0.325\linewidth]{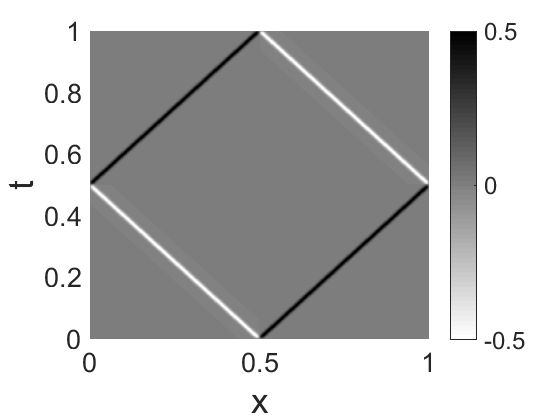}
	\caption{Linear wave equation: snapshots of the full-order solution for the density (left) and the velocity (right).}
	\label{fig:linearWave_FOM_snapshots2}
\end{figure}

\begin{figure}[htb]
	\centering
	\includegraphics [width=0.325\linewidth]{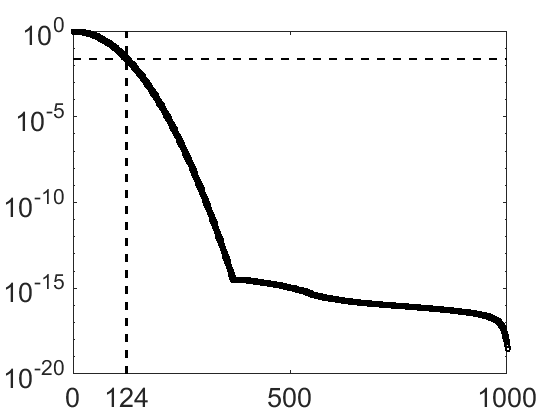}
	\caption{Linear wave equation: singular value decay of the snapshot matrix.}
	\label{fig:svDecay}
\end{figure}

As a model problem for a partial differential equation whose solution features multiple transport velocities we consider the linear acoustic wave equation
\begin{equation}
\begin{aligned}
\partial_t \rho  +  \rhoref   \partial_x u  &=& 0, \\
\partial_t  u   +  c^2/\rhoref \partial_x \rho  &=& 0 ,
\end{aligned}
\label{eq:linearWave}
\end{equation}
on a one-dimensional spatial domain $\Omega= (0,1) $ with periodic boundary conditions.
Here, $u$ is the velocity, $\rho$ the density, $\rhoref$ a reference density, and $c$ the speed of sound.
The analytic solution of \eqref{eq:linearWave} can be expressed as
\begin{equation}
\label{eq:analyticSolution}
\begin{bmatrix}
\rho\left(x,t\right)\\
u\left(x,t\right)
\end{bmatrix}
=
q_{-}\left(x+ct\right)
\begin{bmatrix}
\rhoref\\
-c
\end{bmatrix}
+q_{+}\left(x-ct\right)
\begin{bmatrix}
\rhoref\\
c
\end{bmatrix},
\end{equation}
where $q_{-}$ and $q_{+}$ are the Riemann invariants which are uniquely determined by the initial conditions. In the following, we use $\rhoref=1$ and $c=1$ and we consider the initial conditions
\begin{equation*}
	\rho(x,0)=\rho_0(x)=\exp\left(-\left(\frac{x-0.5}{0.01}\right)^2\right),\qquad u(x,0)=u_0(x)\equiv 0,
\end{equation*}
which represent a pressure pulse with large gradients.
Due to these large gradients, the analytic solution, see Fig.\ \ref{fig:linearWave_FOM_snapshots2}, is very hard to approximate by a classical POD approach, since the singular values of the snapshot matrix, that is obtained by sampling the analytic solution, are decaying very slowly, cf.\ Fig.\ \ref{fig:svDecay}.
To demonstrate the difficulties that POD has for this problem consider the \emph{relative approximation error}
\begin{equation}
	\label{eq:errorMeasure}
	\left.\left(\sum_{j=1}^n \left\lVert \left[X\right]_j-\left[\tilde X\right]_j\right\rVert^2\right)\right/\left(\sum_{j=1}^n \left\lVert \left[X\right]_j\right\rVert^2\right)
\end{equation}
of an approximation $\tilde X$ of the snapshot matrix $X$, $\left\lVert\cdot\right\rVert$ being the Euclidean norm.

In this model problem, to obtain a relative error of less than $1\%$, the POD needs $124$ modes (cf.\ dashed lines in Fig.\ \ref{fig:svDecay}) although the analytic solution is simply represented by the sum of two shifted functions.
Indeed, the analytic solution \eqref{eq:analyticSolution} can be formulated within the more general representation \eqref{eq:shiftedSnapshotDecomposition} with only two modes and
\begin{alignat*}{2}
	&\nShifts=2,\quad r_1=r_2=1,\quad &&\shift^{1}(t)=-\shift^{2}(t)=t,\quad\amplitude^{1}_{1}(t)=\amplitude^{2}_{1}(t)\equiv 0.5,\\
	&\modes^{1}_{1}(x)=\rho_0(x)
	\begin{bmatrix}
		1 & -1
	\end{bmatrix}
	^T,\quad &&\modes^{2}_{1}(x)=\rho_0(x)
	\begin{bmatrix}
		1 &1
	\end{bmatrix}
	^T.
\end{alignat*}
However, the question arises how to compute such a decomposition when only snapshot data are available.
In this case the POD is optimal in the sense that it minimizes the residual, i.\,e., it solves the optimization problem
\begin{equation}
\label{eq:PODoptimization}
\min_{\modes,\amplitude} \int\limits_0^T\int\limits_\Omega \left(q\left(x,t\right)-\sum_{\modeIdx=1}^r\amplitude_\modeIdx(t)\modes_\modeIdx(x)\right)^2 \mathrm{d}x \mathrm{d}t\quad \mbox{s.\ t.} \left\langle\modes_i(x),\modes_j(x)\right\rangle_{L^2(\Omega)}=\delta_{ij}
\end{equation}
for $i,j=1,\ldots,r,$ where $\delta$ denotes the Kronecker delta. In this way the modes $\modes_j$
form an orthonormal basis with respect to the $L^2$ inner product in $\Omega$.

To extend this optimality of \eqref{eq:PODoptimization} to the more general decomposition \eqref{eq:shiftedSnapshotDecomposition}, we consider the optimization problem
\begin{equation}
\label{eq:sPODoptimization_fixedShifts}
\min_{\modes,\amplitude}\int\limits_0^T\int\limits_\Omega \left(q\left(x,t\right)-\sum\limits_{\ell=1}^{\nShifts} \Tperinf\left(\shift^{\ell}(t)\right) \sum\limits_{\modeIdx=1}^{r^\ell} \amplitude^{\ell}_{\modeIdx}(t)\modes^{\ell}_{\modeIdx}(x)\right)^2 \mathrm{d}x \mathrm{d}t,
\end{equation}
where for the moment we assume that the shift frames $\shift$ are available  or can be approximated before the optimization for the modes $\modes$ and their time amplitudes $\amplitude$ is carried out.
Methods to estimate these shifts based on given snapshot data have been discussed in \cite{ReiSSM17_ppt}.

A crucial difference between \eqref{eq:sPODoptimization_fixedShifts} and \eqref{eq:PODoptimization} is that \eqref{eq:sPODoptimization_fixedShifts} is an unconstrained optimization problem without the orthonormality restriction for the modes $\modes_j$. The reason why we drop this orthonormality requirement is that in a decomposition of the form \eqref{eq:shiftedSnapshotDecomposition} even linearly dependent modes may lead to optimal approximations.

To illustrate the necessity to allow linearly dependent modes, consider again the linear wave equation but this time only the density, \ie, take $q(x,t)=\rho(x,t)$. In this case a solution of the optimization problem \eqref{eq:sPODoptimization_fixedShifts} is obtained with
\begin{alignat*}{2}
	&\nShifts=2,\quad r_1=r_2=1,\quad &&\shift^{1}(t)=-\shift^{2}(t)=t,\\
	&\amplitude^{1}_{1}(t)=\amplitude^{2}_{1}(t)\equiv 0.5, \quad &&\modes^{1}_{1}(x)=\modes^{2}_{1}(x)=\rho_0(x),
\end{alignat*}
which means that there exists an optimal approximation with linearly dependent modes $\modes^{1}_{1}=\modes^{2}_{1}$.

For this reason, we prefer not to enforce orthogonality of the modes in \eqref{eq:sPODoptimization_fixedShifts} as in the POD optimization problem \eqref{eq:PODoptimization}, at least when there is more than one transport velocity.

\section{Residual Minimization\label{sec:optimization}}

In this section we discuss the optimization problem \eqref{eq:sPODoptimization_fixedShifts} with a general linear shift operator $\Tinf$, \ie, we consider
\begin{equation}
	\label{eq:generalTransformedModesApproximation}
	\min_{\modes,\amplitude}\int\limits_0^T\int\limits_\Omega \left(q\left(x,t\right)-\sum\limits_{\ell=1}^{\nShifts} \Tinf\left(\shift^{\ell}(t)\right) \sum\limits_{\modeIdx=1}^{r^\ell}
	\amplitude^{\ell}_{\modeIdx}(t)\modes^{\ell}_{\modeIdx}(x)\right)^2 \mathrm{d}x \mathrm{d}t.
\end{equation}
The solution of the POD optimization problem \eqref{eq:PODoptimization} can be obtained by
solving an operator eigenvalue problem,  which in the discrete setting corresponds to computing an SVD.
Since in the setting of \eqref{eq:generalTransformedModesApproximation}, the modes may be linearly dependent, we have to solve a nonlinear optimization problem instead.
To this end, we apply numerical optimization techniques on the discrete level but, prior to that, we analyze some properties of \eqref{eq:generalTransformedModesApproximation}.

First, it should be noted that the solution is in general not unique. This can be seen by taking for instance the simple case where $\Tinf=\Tperinf$ and $q(x,t)=q_1(x+t)+q_2(x-t)+\cos(t)q_3(x)$ with some arbitrary functions $q_i$ for $i=1,\ldots,3$.
Then, a solution of \eqref{eq:generalTransformedModesApproximation} is given by
\begin{alignat*}{2}
	&\nShifts=3,\quad r_1=r_2=r_3=1,\quad &&\shift^{1}(t)=-\shift^{2}(t)=t,\quad \shift^3(t)\equiv 0,\\
	&\amplitude^{1}_{1}(t)=\amplitude^{2}_{1}(t)\equiv 1, \quad &&\amplitude^{3}_{1}(t)=\cos(t),\quad \modes^{i}_{1}(x)=q_i(x),\quad\mbox{for }i=1,\ldots,3.
\end{alignat*}
On the other hand, by making use of the trigonometric identities $\sin(x\pm t)=\sin(x)\cos(t)\pm \cos(x)\sin(t)$, another solution is
\begin{align*}
	&\nShifts=3,\quad r_1=r_2=r_3=1,\quad \shift^{1}(t)=-\shift^{2}(t)=t,\quad \shift^3(t)\equiv 0\\
	&\amplitude^{1}_{1}(t)=\amplitude^{2}_{1}(t)\equiv 1, \quad\amplitude^{3}_{1}(t)= \cos(t),\\
	&\modes^{i}_{1}(x)=q_i(x)+\sin(x),\quad\mbox{for }i=1,2,\quad \modes^{3}_{1}(x) = q_3(x)-2\sin(x).
\end{align*}
Both these solutions are optimal, since the cost functional is zero.

As discussed in Sec.\ \ref{sec:intro}, many of the currently discussed model reduction approaches for transport-dominated phenomena consider the case of only one transport velocity ($\nShifts=1$) and periodic boundary conditions. In this special case the cost functional takes the form
\begin{equation}
	\label{eq:coMovingFrameApproximation}
	\min_{\modes,\amplitude}\int\limits_0^T\int\limits_\Omega \left(q\left(x,t\right)-\Tperinf\left(\shift(t)\right) \sum\limits_{\modeIdx=1}^{r}
	\amplitude_{\modeIdx}(t)\modes_{\modeIdx}(x)\right)^2 \mathrm{d}x \mathrm{d}t
\end{equation}
and one can enforce the modes to form an orthonormal basis, since orthogonality is preserved under the action of the periodic shift operator, \ie,
\begin{equation}
	\label{eq:orthogonalityConstraints}
	\left\langle \Tperinf\left(\shift(t)\right)\modes_i(x),\Tperinf\left(\shift(t)\right)\modes_j(x)\right
\rangle_{L^2(\Omega)}=\left\langle\modes_i(x),\modes_j(x)\right\rangle_{L^2(\Omega)}=\delta_{ij}
\end{equation}
for $i,j=1,\ldots,r$.
This follows, since $\Tperinf(\cdot)$ is a unitary operator, cf. \cite{Coh13}.

Since the adjoint operator of $\Tperinf(\shift)$ is given by $\Tperinf^\ast(\shift)=\Tperinf(-\shift)$, the optimization problem \eqref{eq:coMovingFrameApproximation} associated with the constraints \eqref{eq:orthogonalityConstraints} is equivalent to
\begin{equation*}
	\min_{\modes,\amplitude}\int\limits_0^T\int\limits_\Omega \left(\Tperinf\left(-\shift(t)\right)q\left(x,t\right)-\sum\limits_{\modeIdx=1}^{r}
	\amplitude_{\modeIdx}(t)\modes_{\modeIdx}(x)\right)^2 \mathrm{d}x \mathrm{d}t,\quad \mbox{s.\,t.}\,\mbox{\eqref{eq:orthogonalityConstraints}}.
\end{equation*}
Thus in this special case, the optimization problem leads to a POD of the transformed function $\Tperinf\left(-\shift(t)\right)q\left(x,t\right)$, which has been used, \eg, in \cite{CagMS16_ppt}.

In the general case of more than one transport velocity ($\nShifts>1$), we have to solve the optimization problem \eqref{eq:generalTransformedModesApproximation} numerically. 
Carrying out a discretization, we have to solve the optimization problem
\begin{equation}
	\label{eq:generalTransformedModesApproximation_discrete}
	\min_{\dmodes,\damplitude}\underbrace{\sum_{j=1}^n\left\lVert \left[X\right]_j-\sum\limits_{\ell=1}^{\nShifts} \T\left(\dshift^{\ell}_{j}\right) \sum\limits_{\modeIdx=1}^{r^\ell}
	\damplitude^{\ell}_{\modeIdx,j}\dmodes^{\ell}_{\modeIdx}\right\rVert^2}_{=\vcentcolon \J}.
\end{equation}
where $n$ is the number of snapshots.
Introducing the notation
\begin{align*}
\nonumber 
	 &\damplitude_j\vcentcolon=[\damplitude^{1}_{1,j}\quad\ldots\quad\damplitude^{1}_{r_1,j}\quad
\damplitude^{2}_{1,j}\quad\ldots\quad\damplitude^{\nShifts}_{r_\nShifts,j}]^T,\\
	&K_j\vcentcolon=
	\begin{bmatrix}
		T\left(\dshift^{1}_{j}\right)\dmodes^{1}_{1} & \ldots &\T\left(\dshift^{1}_{j}\right)\dmodes^{1}_{r_1} & \T\left(\dshift^{2}_{j}\right)\dmodes^{2}_{1} & \ldots & \T\left(\dshift^{\nShifts}_{j}\right)\dmodes^{\nShifts}_{r_\nShifts}
	\end{bmatrix}	
	,
\end{align*}
the cost functional in \eqref{eq:generalTransformedModesApproximation_discrete} can be expressed as the least squares problem
\begin{equation}
\label{eq:shortCostFunctional}
\J = \sum\limits_{j=1}^n\left\lVert \left[X\right]_{j}-K_j\damplitude_j\right\rVert_2^2.
\end{equation}
Considering the dependency of $J$ with respect to the amplitudes $\damplitude_j$ for fixed modes $\dmodes$, the necessary optimality condition is given by
\[
	\nabla_{\damplitude_j}\J = -2K_j^T\left(\left[X\right]_{j}-K_j\damplitude_j\right)=0,
\]
or equivalently
\begin{equation}
	\label{eq:alphaCondition}
K_j^TK_j\damplitude_j = K_j^T\left[X\right]_{j}
\end{equation}
for $j=1,\ldots,n$.
The general solution of \eqref{eq:alphaCondition} is given by
\begin{equation}
	\label{eq:solutionOfLS}
	\damplitude_j = V_{j,1}\Sigma_{j,1}^{-1}U_{j,1}^T\left[X\right]_{j}+V_{j,2}\beta_j,
\end{equation}
where $\beta_j$ is an arbitrary vector of suitable dimension, and the matrices $V_{j,1}$, $\Sigma_{j,1}$, $U_{j,1}$, and $V_{j,2}$ are defined via the SVD of $K_j$
\begin{equation*}
	K_j =
	\begin{bmatrix}
		U_{j,1} & U_{j,2}
	\end{bmatrix}
	\begin{bmatrix}
		\Sigma_{j,1} & 0\\
		0 & 0
	\end{bmatrix}
	\begin{bmatrix}
		V_{j,1}^T \\ V_{j,2}^T
	\end{bmatrix},
\end{equation*}
where $\Sigma_{j,1}$ contains the non-zero singular values of $K_j$ \cite{GanGK14}.
If the shifted modes are linearly independent at a time point $t$, then $V_{j,2}$ is void and the solution \eqref{eq:solutionOfLS} is unique, otherwise
\eqref{eq:alphaCondition} has infinitely many solutions.

Substituting \eqref{eq:solutionOfLS} into \eqref{eq:shortCostFunctional}, the cost functional takes the form
\begin{equation*}
	J = \sum\limits_{j=1}^n\left\lVert \left[X\right]_{j}-U_{j,1}U_{j,1}^T\left[X\right]_{j}\right\rVert_2^2,
\end{equation*}
which only depends on the modes $\dmodes$ hidden in the matrices $U_{j,1}$.
Simple calculations show that minimizing $\J$ is equivalent to the minimization problem
\begin{equation}
	\label{eq:finalOptimizationProblem}
	\min_{\dmodes}
\tilde \J = -\sum\limits_{j=1}^n\left\lVert U_{j,1}^T\left[X\right]_{j}\right\rVert_2^2.
\end{equation}
The gradient of $\tilde J$ with respect to a mode $\dmodes^{\ell}_{k}$ is given by
\begin{equation*}
	\label{eq:gradientOfJ}
	\nabla_{\dmodes^{\ell}_{k}}\tilde \J = \sum\limits_{j=1}^n \damplitude^{\ell}_{k,j} \T\left(\dshift^{\ell}_{j}\right)^T\left(I_{m}-U_{j,1}U_{j,1}^T\right)\left[X\right]_{j}.
\end{equation*}
An algorithm to compute an optimal solution is presented in Sec.\ \ref{sec:algorithm}.

\section{Algorithm and Implementation\label{sec:algorithm}}

Since it is a priori unclear how many modes are necessary to achieve a certain error tolerance, we propose to solve the optimization problem \eqref{eq:finalOptimizationProblem} starting with a small number of modes and iteratively adding modes in a greedy fashion, cf. Algorithm \ref{alg:sPOD}.
\begin{algorithm}[htb]
	\caption{sPOD algorithm based on residual minimization}
	\label{alg:sPOD}
	\begin{algorithmic}[1]
		\Statex \textbf{Input:} snapshot matrix $X$; initial mode numbers $\rVec^0$; shifts $\dshift^{\ell}_{j}$ for $j=1,\ldots,n$ and $\ell=1,\ldots,\nShifts$; routine for the calculation of $\T(\cdot)$; error
		tolerance \textit{tol}; maximum iteration number $\iter_{\mathrm{max}}$
		\Statex \textbf{Output:} modes $\dmodes^{\ell}_{\modeIdx}$; amplitudes $\damplitude^{\ell}_{\modeIdx,j}$ for $j=1,\ldots,n$, $\ell=1,\ldots,\nShifts$, and $k=1,\ldots,r^{\ell}$
		\State Solve \eqref{eq:finalOptimizationProblem} with mode numbers $\rVec^0$ for the modes $\dmodes$ \label{line:firstOptimization}
		\State Compute the amplitudes $\damplitude$ from \eqref{eq:solutionOfLS} \label{line:alpha}
		\State Reconstruct $\tilde X$ as in \eqref{eq:shiftedSnapshotDecomposition_discrete} and compute the error as in \eqref{eq:errorMeasure} \label{line:error}
		\State $\iter=0$
		\While{$(\text{error}>\text{tol})\ \mbox{\rm and } (\iter<\iter_{\mathrm{max}})$}\label{line:while}
			\State $\iter \gets \iter+1$
			\For{$i=1:\nShifts$} \label{line:for}
				\State Solve \eqref{eq:finalOptimizationProblem} with mode numbers $\rVec^{\iter-1}+e_i$ for the modes $\dmodes$ \label{line:optimizationInForLoop}
				\State Compute the amplitudes $\damplitude$, $\tilde X$, and the error as in lines \ref{line:alpha} and \ref{line:error}
				\State $\text{tempError}(i)\gets\text{error}$
			\EndFor \label{line:endFor}
			\State Find the index $q$ for which $\text{tempError}$ is minimal
			\State $\text{error}\gets\text{tempError}(q) $
			\State $\rVec^\iter\gets \rVec^{\iter-1}+e_q$
		\EndWhile
	\end{algorithmic}
\end{algorithm}
To initiate the algorithm we choose a vector $\rVec^0\in\mathbb{N}^{\nShifts}$ containing the initial mode numbers for each velocity frame and prescribed shifts $\dshift^{\ell}_{j}$ for each velocity frame and discrete time step.
The algorithm starts with computing a mode decomposition with mode numbers $\rVec^0$.
For this, the optimization problem \eqref{eq:finalOptimizationProblem} is solved with a nonlinear optimization solver of choice, \eg, Newton's method or quasi-Newton methods, see \eg~\cite{NocW06}.
Since the optimization problem scales with the full dimension, 
we recommend an inexact Newton method or a limited-memory quasi-Newton method which is more efficient \cite{NocW06}.
Motivated by the case with one velocity frame and a periodic shift operator discussed in Sec.\ \ref{sec:optimization}, we choose the first $[\rVec^0]_\ell$ singular vectors of the transformed snapshot matrix
\begin{equation*}
	\begin{bmatrix}
		\T\left(-\dshift^{\ell}_{1}\right)[X]_1 & \cdots & \T\left(-\dshift^{\ell}_{n}\right)[X]_n
	\end{bmatrix}
\end{equation*}
as starting values for the modes of the $\ell$th velocity frame.
Following this, in line \ref{line:while} the relative error is compared with the tolerance and if the tolerance is not achieved, then the algorithm continues by adding modes in a greedy manner.
More precisely, in the \textbf{for} loop in lines \ref{line:for}-\ref{line:endFor}, we add one mode to each frame at a time, solve the optimization problem \eqref{eq:finalOptimizationProblem}, construct $\tilde X$, and compute the error.
Subsequently, the errors corresponding to the different mode number vectors $\rVec^{\iter-1}+e_i$ are compared, where $e_i\in\mathbb{R}^\nShifts$ denotes the $i$th unit vector, and only that mode is kept which results in the smallest error.
This \textbf{while} loop continues until the error is below the tolerance or the maximum iteration number is reached.

The major computational cost of Algorithm~\ref{alg:sPOD} arises from the solution of the optimization problems in lines \ref{line:firstOptimization} and \ref{line:optimizationInForLoop} and depends on the chosen solver.
The computation time can be decreased significantly by performing the \textbf{for} loop in lines \ref{line:for}-\ref{line:endFor} in parallel.
Another opportunity for a speedup is to use multigrid methods for the optimization, see \eg{} \cite{Nas00}.

Most parts of Sec.\ \ref{sec:optimization}, as well as Algorithm \ref{alg:sPOD} are valid for general matrix functions $\T$ which do not necessarily have to be associated with a shift operation. Thus, the use of matrix functions which simulate other effects like rotation or dilation is possible, however, this topic is not within the scope of this paper.

Instead, in Sec.\ \ref{sec:DDT} we use a shift operator with constant extrapolation, \ie,
\begin{equation*}
	\Tcinf\left(\shift(t)\right)f\left(x,t\right)\vcentcolon=\left\{
	\begin{array}{lll}
		f\left(x-\Delta(t),t\right) & \text{for } 0\leq x-\Delta(t) \leq L,\\
		f\left(0\right) & \text{for } x-\Delta(t)<0,\\
		f\left(L\right) & \text{for } x-\Delta(t)>L.
	\end{array}
    \right.
\end{equation*}
For the discrete analogue $\Tc$ on a uniform grid with mesh width $h$, we distinguish between two cases: If the shift is a multiple of $h$, then $\Tc(\cdot)$ is defined as
\begin{equation}
	\label{eq:shiftOperatorWithConstantExtrapolation}
	\Tc\left(kh\right) = \left[
	\begin{array}{cccc|c}
		1 & 0 & \ldots & 0 & 0 \\
		\hline
		 & & & & 0\\
		 & & & & \vdots\\
		 \multicolumn{4}{c|}{\text{\smash{\raisebox{1.5\normalbaselineskip}{$I_{m-1}$}}}} & 0		 
	\end{array}
	\right]^k,\quad
	\Tc\left(-kh\right) = \left[
	\begin{array}{c|cccc}
		0 & & & & \\
		\vdots & & & & \\
		0 & \multicolumn{4}{|c}{\text{\smash{\raisebox{1.5\normalbaselineskip}{$I_{m-1}$}}}}\\\hline
		0 & 0 & \cdots & 0 & 1		
	\end{array}
	\right]^k
\end{equation}
with $k\in\mathbb{N}$.
If the shift is not a multiple of $h$ we use an interpolation scheme, \ie, for instance, a linear interpolation like $\Tc(0.5h)=0.5(\Tc(0)+\Tc(h))$.
Similarly, a shift matrix function for the periodic case has been introduced in \cite{RimML17_ppt}. 

\section{Test Case: Pulsed Detonation Combuster\label{sec:DDT}}
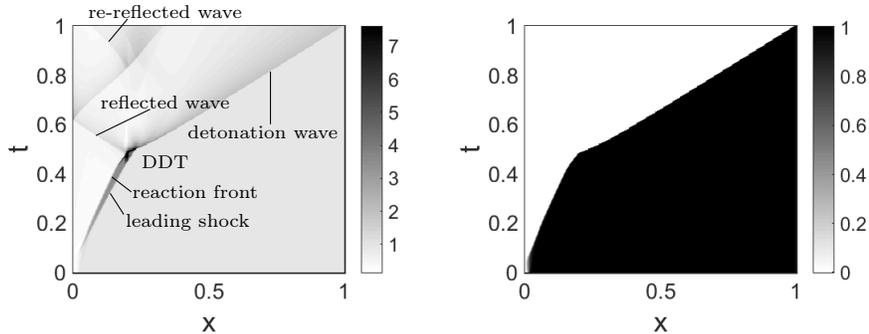
\begin{figure}[htb]
	\centering
	\begin{tikzpicture}
		\node[inner sep=0pt] (rho_FOM) at (-3,0){\includegraphics[width=0.49\linewidth]{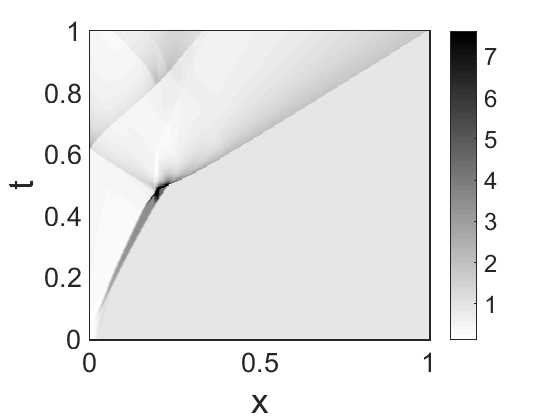}};
		\node[inner sep=0pt] (Y_FOM) at (3,0){\includegraphics[width=0.49\linewidth]{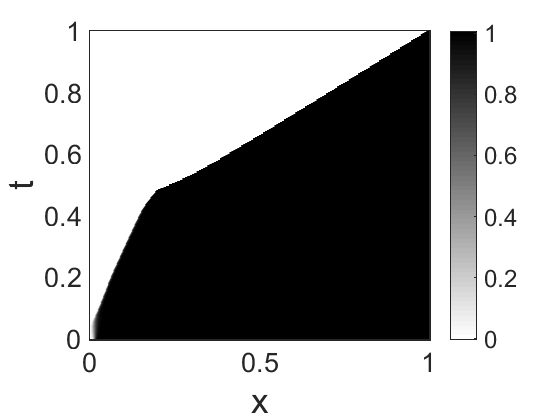}};
		\node[inner sep=0pt] (DDT) at (-3.8,0.1) {\scriptsize{DDT}};
		\node[inner sep=0pt] (LS) at (-3.5,-0.7) {\scriptsize{leading shock}};
		\node[inner sep=0pt] (RF) at (-3.4,-0.3) {\scriptsize{reaction front}};
		\node[inner sep=0pt] (RW) at (-3.8,0.9) {\scriptsize{reflected wave}};
		\node[inner sep=0pt] (RRW) at (-3.8,2.1) {\scriptsize{re-reflected wave}};
		\node[inner sep=0pt] (DW) at (-2.5,0.5) {\scriptsize{detonation wave}};
		\draw (-4.55,1.7) -- (-4.3,2);
		\draw (-4.72,0.45) -- (-4,0.8);
		\draw (-2.4,1.3) -- (-2.4,0.6);
		\draw (-4.51,-0.1) -- (-4.25,-0.3);
		\draw (-4.53,-0.32) -- (-4.35,-0.6);
	\end{tikzpicture}
	\caption{PDC: snapshots of the full-order solution for density (left) and species (right).}
	\label{fig:FOM_snapshots}
\end{figure}

As a realistic test example, we consider density, velocity, pressure, and effective species snapshot data of a \emph{Pulsed Detonation Combuster (PDC)} with a shock-focusing geometry.
The data is based on a simulation of the reactive, compressible Navier-Stokes equations where physical parameters have been adjusted by a data assimilation, see \cite{GraLRPSM17}.
The density and species snapshots are depicted in Fig.\ \ref{fig:FOM_snapshots}.
In the snapshots of the species we observe a reaction front propagating through the domain.
The density snapshots show initially two transports, the reaction front and a leading shock, slightly diverging before they converge again and interact.
This \emph{deflagration to detonation transition (DDT)} is caused by a nozzle at around $x=0.2$, cf. \cite{GraLRPSM17}.
Following this, the reaction front and the leading shock continue as a detonation wave moving to the right.
At the same time, a reflected wave is moving to the left before being reflected at the boundary.
When it reaches the nozzle again, another partial reflection is visible.
The velocity and pressure snapshots look similar.

\begin{figure}[b]
	\centering
	\includegraphics [width=0.325\linewidth]{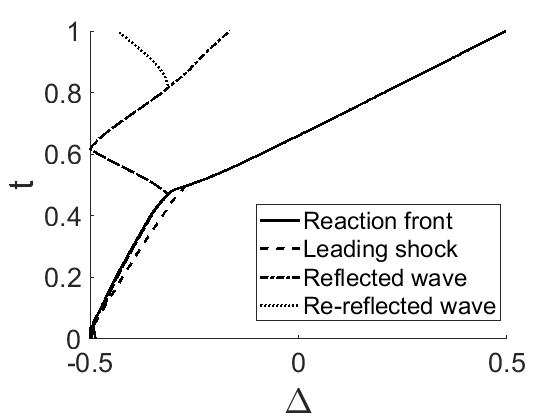}
	\caption{PDC: tracked shifts for the different transports.}
	\label{fig:cs}
\end{figure}

Before we apply Algorithm \ref{alg:sPOD} we need to find good candidates for the shifts corresponding to the transports of the system.
Here, we focus on the four most dominant transports: the reaction front, the leading shock, the reflected wave, and the partial reflection at the nozzle which is referred to as \textit{re-reflected wave} in the following.
We track these transports based on the snapshot data without any a priori knowledge of their velocities.
The reaction front is the easiest to detect, since it is clearly visible in the species snapshots as a large gradient.
To track it, we determine the location of the maximum in each column of the difference matrix whose $j$th column is defined as the difference between the $j+1$st and $j$th column of the species snapshot matrix.
The resulting tracked shift is depicted in Fig.\ \ref{fig:cs}, solid line.
Here, negative shift values occur since the reaction front is shifted such that it is centered in the middle of the computational domain.

The tracking of the other transports works similarly, but is a little more elaborate since we need to distinguish them from each other.
To this end, we restrain the region of the computational domain where the location of the maximum slope is computed. This subregion depends on both the considered transport and time interval.
In our tests, this decomposition in subregions has been done manually based on the velocity snapshots. The corresponding tracked shifts are depicted in Fig.\ \ref{fig:cs}.
In addition, we also add a frame with zero velocity to account for the structures that we cannot capture well by the other velocity frames.

We apply Algorithm \ref{alg:sPOD} with a shift operator with constant extrapolation as in \eqref{eq:shiftOperatorWithConstantExtrapolation} with Lagrange polynomials of degree three for the interpolation.
In addition, we specify $\mathrm{tol}=0.01$, $\iter_{\mathrm{max}}=n$, and $\rVec^0=[1\;1\;1\;1\;0]$, \ie, one mode for each of the non-zero velocity frames.
The nonlinear optimization problem is solved using the MATLAB package HANSO which is based on a limited-memory BFGS \cite{LewO13}.
Moreover, to avoid parasitic structures in the approximation of the species, we force those parts of the modes which correspond to the species and to other transports than the reaction front, to be zero.

In this test case we  have to deal with data of physical variables with highly different scales.
To avoid that the approximation of the physical variable with the highest scale becomes dominant we scale the snapshots such that the snapshot matrices of the different physical variables have the same Frobenius norm.
We build the snapshot matrix $X$ for Algorithm \ref{alg:sPOD} by concatenating the scaled snapshot matrices of the different physical variables.

\begin{figure}[t]
	\centering
	\includegraphics [width=0.325\linewidth]{rho_FOM}
	\includegraphics [width=0.325\linewidth]{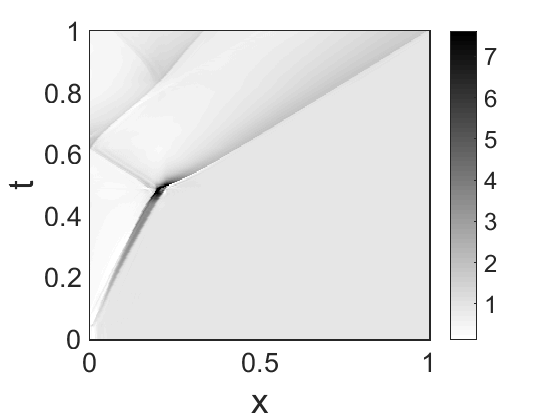}
	\includegraphics [width=0.325\linewidth]{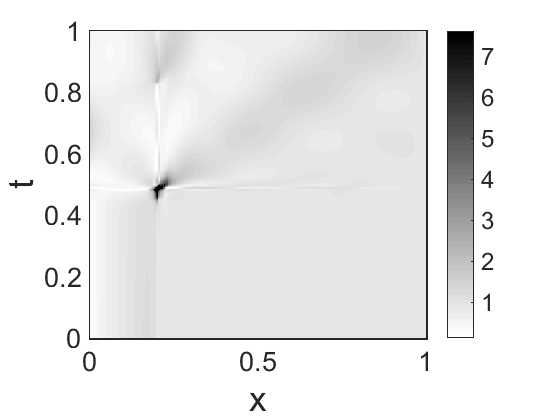}
	\includegraphics [width=0.325\linewidth]{Y_FOM}
	\includegraphics [width=0.325\linewidth]{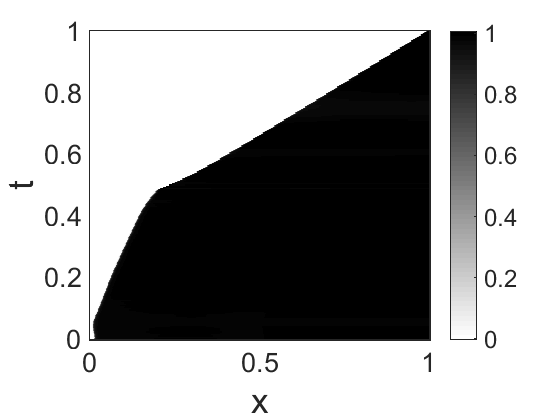}
	\includegraphics [width=0.325\linewidth]{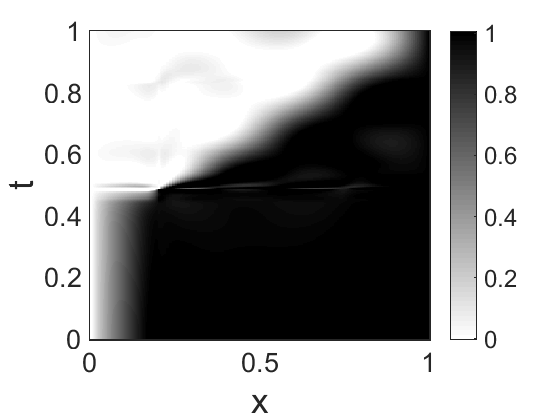}
	\caption{PDC: comparison between full-order solution (left column), sPOD approximation with $7$ modes (middle column), and POD approximation with $7$ modes (right column). The top
	row shows the results for the density, the bottom row for the species.}
	\label{fig:FOMvsSPODvsPOD}
\end{figure}

Algorithm \ref{alg:sPOD} terminates after $3$ iterations in the {\bf while} loop with an error of $0.71\%$ and $\rVec^3=[3\;1\;1\;1\;1]$, \ie, two modes have been added to the reaction front and one mode to the zero velocity frame. This means that we meet the error tolerance with $7$ modes in total. The sPOD approximation for the density and the species is depicted in Fig.\ \ref{fig:FOMvsSPODvsPOD}, middle column.
Although some deviations to the full-order solution are visible, the sPOD captures the dynamics well and the dominant transports are clearly distinct.
This becomes even more striking when comparing it to the POD with the same number of modes which is plotted in Fig.\ \ref{fig:FOMvsSPODvsPOD}, right.
As is common in the POD literature, we first subtracted the mean value of each row of the snapshot matrix to center the data around the origin, cf. \cite{Cha00}.
The POD approximation of the density features a high peak in the region of the DDT while the other structures are hardly recognizable.
For the species, the reaction front is at least indicated, but blurred, and further distortions are visible especially near the DDT.
To obtain a POD approximation of the same accuracy as the sPOD with $7$ modes, $73$ POD modes are needed for this example.

\begin{figure}[t]
	\centering
	\includegraphics [width=0.325\linewidth]{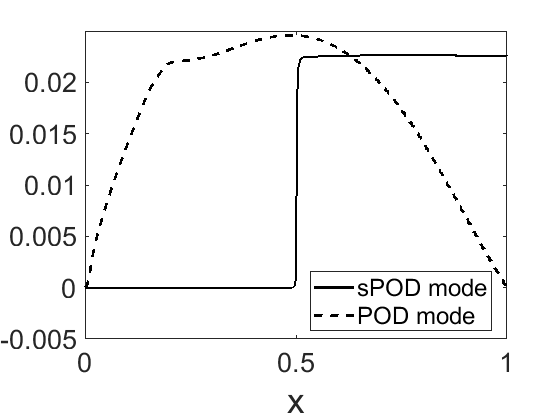}
	\caption{PDC: comparison of first POD mode and first sPOD mode for the species.}
	\label{fig:Ymodes}
\end{figure}

Another advantage of the sPOD becomes clear when looking at the POD and sPOD modes.
In Fig.\ \ref{fig:Ymodes} the first sPOD mode for the species in the reaction front frame is depicted and compared to the first POD mode.
While the sPOD mode clearly reveals the reaction front as a jump in the middle, the POD is rather smooth and does not show any structure resembling a reaction front.

\begin{figure}
	\centering
	\includegraphics [width=0.325\linewidth]{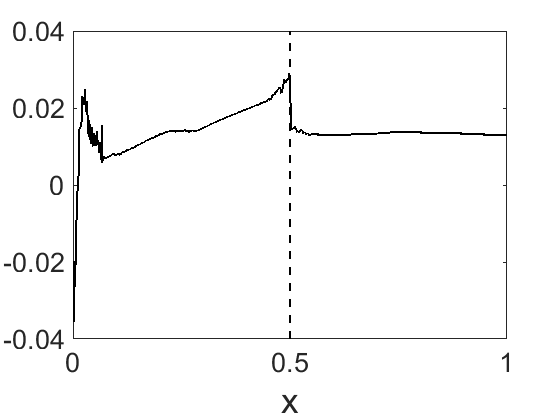}
	\includegraphics [width=0.325\linewidth]{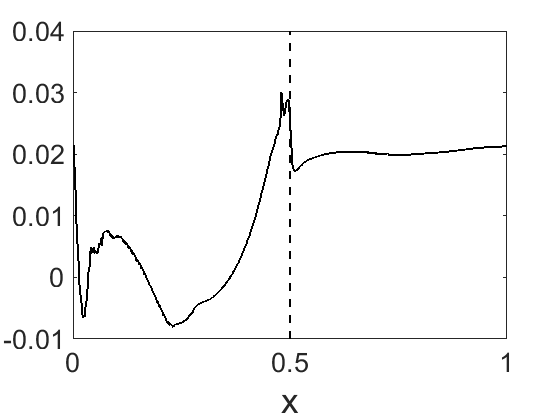}
	\includegraphics [width=0.325\linewidth]{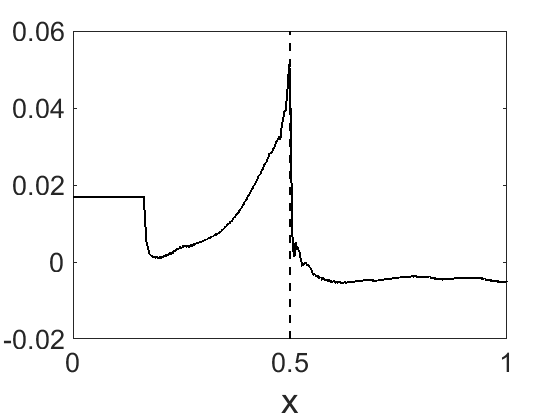}
	\includegraphics [width=0.325\linewidth]{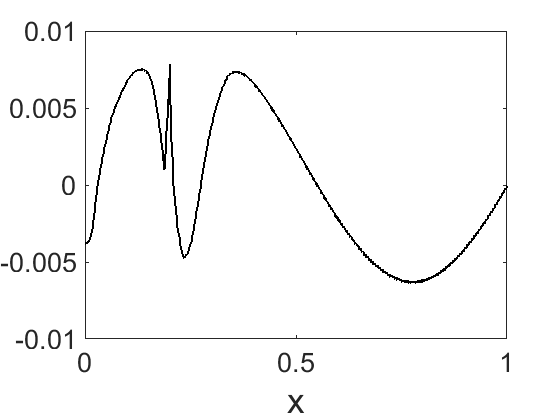}
	\includegraphics [width=0.325\linewidth]{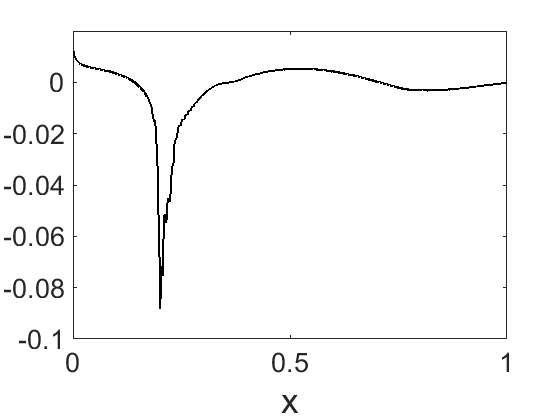}
	\includegraphics [width=0.325\linewidth]{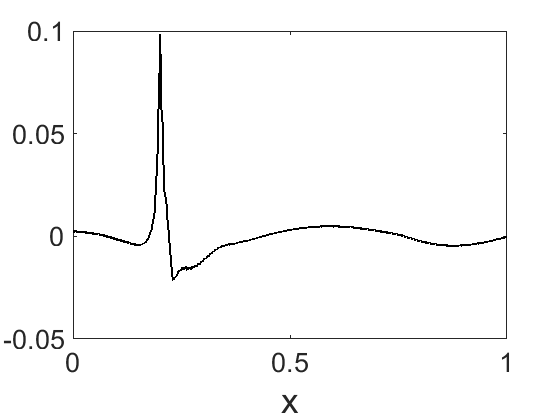}
	\caption{PDC: comparison of first sPOD modes (top row) for the reaction front, the leading shock, and the reflected wave (from left to right) and first three POD
	modes (bottom row, from left to right) for the density.}
	\label{fig:rhoModes}
\end{figure}

In Fig.\ \ref{fig:rhoModes} the first sPOD mode for the density is depicted for the reaction front, the leading shock, and the reflected wave and compared to the first three POD modes.
The latter ones mainly focus on the DDT which agrees with Fig.\ \ref{fig:FOMvsSPODvsPOD}, top right, while the shapes of the moving fronts are not captured.
The sPOD modes are not as clear as in Fig.\ \ref{fig:Ymodes} but still each of them features a clear front profile in the middle (marked by dashed lines) which reflects the physics properly.
However, especially at the left boundary they differ strongly: The mode for the reflected wave, top right in Fig.\ \ref{fig:rhoModes}, reveals a flat profile at the left boundary.
This is due to the fact that this part of the mode is not used in the sPOD approximation, since the corresponding shift, depicted in Fig.\ \ref{fig:cs}, does not attain values greater than $-0.16$.
The modes for the reaction front and the leading shock reveal some oscillations at the left boundary.
A possible reason for this is the use of the shift operator with constant extrapolation which provides the values at the boundaries of the mode with a disproportionate weight.

\section{Summary and Outlook\label{sec:conclusion}}

We have presented a new algorithm for computing a shifted proper orthogonal decomposition (sPOD) based on a residual minimization applied to snapshot data.
We have applied the algorithm to snapshots  determined from a pulsed detonation combuster (PDC) and compared the results with the standard proper orthogonal decomposition (POD).
The sPOD yields a reasonable approximation of the snapshots with only very few modes.
In contrast, the POD approximation with the same number of modes is blurred and the dynamics is not captured well.
Moreover, the sPOD modes clearly reveal the front profiles of the different transports, whereas the POD is not suitable for identifying structures in this test case.
In comparison to the heuristic sPOD algorithm proposed in \cite{ReiSSM17_ppt}, the new algorithm is based on a residual minimization and hence at least locally optimal.
A drawback of the new algorithm is that it is more expensive than the  POD
and the original sPOD approach of \cite{ReiSSM17_ppt}.
The reason is that a large-scale nonlinear optimization problem has to be solved.

The results of the new sPOD algorithm look promising in terms of both the number of required modes and the physical structures identified by the sPOD modes.
The use of the sPOD modes to obtain a reduced-order model via projection is currently under investigation. Further interesting research directions are a rigorous treatment of non-periodic boundary conditions and an optimization of the shifts together with the modes and the amplitudes.\\

\noindent\textbf{Acknowledgements.} The authors gratefully acknowledge the support by the Deut\-sche Forschungsgemeinschaft (DFG) as part of the collaborative research center SFB 1029 \textit{Substantial efficiency increase in gas turbines through direct use of coupled unsteady combustion and flow dynamics}, project A02 \textit{Development of reduced-order models for pulsed combustion}.

\bibliographystyle{abbrv}
\bibliography{Refs}

\begin{thebibliography}{10}

\bibitem{AbgAC16}
R.~Abgrall, D.~Amsallem, and R.~Crisovan.
\newblock Robust model reduction by {$L^1$}-norm minimization and approximation
  via dictionaries: application to nonlinear hyperbolic problems.
\newblock {\em Adv. Model. Simul. Eng. Sci.}, 3(1), 2016.

\bibitem{Ant05}
A.~C. Antoulas.
\newblock {\em Approximation of large-scale dynamical systems}.
\newblock SIAM, Philadelphia, USA, 2005.

\bibitem{BenGW15}
P.~Benner, S.~Gugercin, and K.~Willcox.
\newblock A survey of projection-based model reduction methods for parametric
  dynamical systems.
\newblock {\em SIAM review}, 57(4):483--531, 2015.

\bibitem{BerHL93}
G.~Berkooz, P.~Holmes, and J.~L. Lumley.
\newblock The proper orthogonal decomposition in the analysis of turbulent
  flows.
\newblock {\em Annu. Rev. Fluid Mech.}, 25:539--575, 1993.

\bibitem{CagMS16_ppt}
N.~Cagniart, Y.~Maday, and B.~Stamm.
\newblock Model order reduction for problems with large convection effects.
\newblock Preprint hal-01395571, 2016.
\newblock Available from \emph{https://hal.archives-ouvertes.fr}.

\bibitem{Cha00}
A.~Chatterjee.
\newblock An introduction to the proper orthogonal decomposition.
\newblock {\em Current Sci.}, 78(7):808--817, 2000.

\bibitem{Coh13}
L.~Cohen.
\newblock {\em The Weyl Operator and its Generalization}.
\newblock Springer Basel, Switzerland, 2013.

\bibitem{GanGK14}
W.~Gander, M.~J. Gander, and F.~Kwok.
\newblock {\em Scientific Computing}.
\newblock Springer International Publishing, Cham, Switzerland, 2014.

\bibitem{GerL14}
J.-F. Gerbeau and D.~Lombardi.
\newblock Approximated {Lax} pairs for the reduced order integration of
  nonlinear evolution equations.
\newblock {\em J. Comput. Phys.}, 265:246--269, 2014.

\bibitem{GraLRPSM17}
J.~A.~T. Gray, M.~Lemke, J.~Reiss, C.~O. Paschereit, J.~Sesterhenn, and J.~P.
  Moeck.
\newblock A compact shock-focusing geometry for detonation initiation:
  {e}xperiments and adjoint-based variational data assimilation.
\newblock {\em Combust. Flame}, 183:144--156, 2017.

\bibitem{HesRS16}
J.~S. Hesthaven, G.~Rozza, and B.~Stamm.
\newblock {\em Certified reduced basis methods for parametrized partial
  differential equations}.
\newblock Springer International Publishing, Cham, Switzerland, 2016.

\bibitem{HolLB96}
P.~Holmes, J.~L. Lumley, and G.~Berkooz.
\newblock {\em Turbulence, Coherent Structures, Dynamical Systems and
  Symmetry}.
\newblock Cambridge University Press, Cambridge, UK, 1996.

\bibitem{IolL14}
A.~Iollo and D.~Lombardi.
\newblock Advection modes by optimal mass transfer.
\newblock {\em Phys. Rev. E}, 89(2):022923, 2014.

\bibitem{KutBBP16}
J.~N. Kutz, S.~L. Brunton, B.~W. Brunton, and J.~L. Proctor.
\newblock {\em Dynamic Mode Decomposition}.
\newblock SIAM, Philadelphia, USA, 2016.

\bibitem{LewO13}
A.~S. Lewis and M.~L. Overton.
\newblock Nonsmooth optimization via quasi-{N}ewton methods.
\newblock {\em Math. Program.}, 141(1--2):135--163, 2013.

\bibitem{MojB17}
R.~Mojgani and M.~Balajewicz.
\newblock Lagrangian basis method for dimensionality reduction of convection
  dominated nonlinear flows.
\newblock Preprint 1701.04343v1, ArXiv, 2017.

\bibitem{Nas00}
S.~G. Nash.
\newblock A multigrid approach to discretized optimization problems.
\newblock {\em Optim. Methods Softw.}, 14(1--2):99--116, 2000.

\bibitem{NocW06}
J.~Nocedal and S.~J. Wright.
\newblock {\em Numerical Optimization}.
\newblock Springer New York, USA, second edition, 2006.

\bibitem{OhlR13}
M.~Ohlberger and S.~Rave.
\newblock Nonlinear reduced basis approximation of parameterized evolution
  equations via the method of freezing.
\newblock {\em C. R. Math. Acad. Sci. Paris}, 351(23--24):901--906, 2013.

\bibitem{ReiSSM17_ppt}
J.~Reiss, P.~Schulze, J.~Sesterhenn, and V.~Mehrmann.
\newblock The shifted proper orthogonal decomposition: {A} mode decomposition
  for multiple transport phenomena.
\newblock Preprint 1512.01985v2, ArXiv, 2017.

\bibitem{RimML17_ppt}
D.~Rim, S.~Moe, and R.~J. LeVeque.
\newblock Transport reversal for model reduction of hyperbolic partial
  differential equations.
\newblock Preprint 1701.07529v1, ArXiv, 2017.

\bibitem{SchS08}
P.~J. Schmid and J.~L. Sesterhenn.
\newblock Dynamic mode decomposition of numerical and experimental data.
\newblock In {\em 61st Annual Meeting of the APS Division of Fluid Dynamics},
  page 208, San Antonio, USA, 2008.

\bibitem{SesS16}
J.~Sesterhenn and A.~Shahirpour.
\newblock A {L}agrangian dynamic mode decomposition.
\newblock Preprint 1603.02539v1, ArXiv, 2016.

\bibitem{Vol01}
S.~Volkwein.
\newblock Optimal control of a phase-field model using proper orthogonal
  decomposition.
\newblock {\em ZAMM Z. Angew. Math. Mech.}, 81(2):83--97, 2001.

\end{thebibliography}

\end{document}